\documentclass[11pt]{report}
\usepackage{amsmath, amscd, amssymb, amsthm}
\usepackage[subnum]{cases}
\usepackage{changepage}
\usepackage{marginnote}
\usepackage{hyperref}
\usepackage{cleveref}
\usepackage[all]{xy}
\usepackage{tabularx}
\usepackage{ltablex}
\usepackage[T1]{fontenc}
\usepackage{setspace}
\usepackage{fancybox}
\usepackage{ragged2e}
\usepackage[margin=1.1in]{geometry}
\hypersetup{colorlinks,linkcolor={black},citecolor={black}}

\begin{document}

\title{$\infty$-Categorical Functional Analysis and $p$-adic Motives}
\author{Xin Tong}
\date{}

\maketitle

\begin{abstract}
We discuss the deep relationship between $\infty$-categorical functional analysis and the anticipated theory of $p$-adic motives. The motivation fundamentally comes from applications essentially in arithmetics from very broad perspectives. The $\infty$-categorical homotopical aspects we considered here come mainly from Bambozzi-Ben-Bassat-Kremnizer, Ben-Bassat-Mukherjee, Bambozzi-Kremnizer, Clausen-Scholze and Kelly-Kremnizer-Mukherjee, based on deep robust foundation of $\infty$-categorical solids and $\infty$-categorical ind-Banach or bornological modules over general Banach rings.
\end{abstract}

\tableofcontents

\chapter{Introduction to Noncommutative Tamagawa Number Conjectures}

\section{The \'Etale Picture}

\noindent We first give the corresponding discussion around the fundamental idea around and introduction to Noncommutative Tamagawa Number Conjecture for motives in the following.
	
\begin{itemize}
\item<1-> 	We follow notations in Burns and Venjakob's paper \cite[Section 4.2, Section 4.3]{BV}, let $M$ be a $K$-Motive, where $K/\mathbb{Q}$ finite extension;

\item<2-> 	$U=\mathrm{Spec}(\mathbb{Z}[1/S])$, where $S$ is a finite set of primes of $\mathbb{Z}$ containing $p$ and $\infty$;
\item<3-> 	$G$ a $p$-adic Lie group with the Iwasawa algebra $\Lambda$, we assume that $G$ factor through the Galois group $G_\mathbb{Q}$, namely we have a $G$-extension $F_\infty/\mathbb{Q}$;
\item<4-> 	Let $\mathbb{T}:=T_p\otimes \Lambda$, the Galois group $G_S$ acts through the quotient and diagonally, where $T_p$ is a $\mathbb{Z}_p$-lattice of $M_p$;
\item<5-> 	Let $\lambda$ be a finite prime of $K$, let $\rho:G\rightarrow \mathrm{GL}_n(\mathcal{O}_\lambda)$ be a continuous representation which is assumed to come from a $K$-motive $N$. We set $M(\rho^*)$ to be $N^*\otimes M$.
\end{itemize}

\begin{itemize}
\item<1-> (\text{Conjecture, Fukaya-Kato 06 \cite[Conjecture 2.3.2]{FK}, Burns-Venjakob \cite[Conjecture 4.1]{BV}}) There is a unique canonical way to associate to the data above an isomorphism in the category $\mathcal{C}_\Lambda$:
\[
\xymatrix@R+0pc@C+0pc{
\zeta_{\Lambda,T}:1_\Lambda \ar[r]\ar[r]\ar[r] & \mathrm{Determinant}_\Lambda (R\Gamma_{c}(U,\mathbb{T}))^{-1}, 
}
\]
such that the base change to $K_\lambda^n$ gives rise to the following isomorphism:
\[
\xymatrix@R+0pc@C+0pc{
1_{K_\lambda} \ar[r]\ar[r]\ar[r] & \triangle_{K}(M(\rho^{*}))_{K_\lambda} \ar[r]\ar[r]\ar[r] & \mathrm{Determinant}_{K_\lambda}(R\Gamma_{c}(U,M(\rho^{*})_{\lambda}))^{-1}, 
}
\]
where the first map is a unique isomorphism such that:
\[
\xymatrix@R+0pc@C+0pc{
1_{\infty} \ar[r]\ar[r]\ar[r] & \triangle_{K}(M({\rho}^{*}))_{\infty} \ar[r]\ar[r]\ar[r] & 1_{\infty}
}
\]
gives rise to the leading term $L^*(M({\rho}^*))$ of the $L$-function of the motive (after considering an embedding of $K$ into $\mathbb{C}$). The second map is the conjectural regulator map predicted by the corresponding standard conjecture. We use here the formulation from Burns and Venjakob in their leading term paper \cite{BV}, but this essentially dates back to Fukaya-Kato \cite[Conjecture 2.3.2]{FK}.

\item<2-> This is basically far-reaching due to many foundational issue, although this is the key ingredients for the construction of $p$-adic L-functions as in Fukaya-Kato's paper \cite{FK}.

\end{itemize}

\noindent{\text{Example: $p$-adic Kubota-Leopoldt and Iwasawa Main Conjecture}}

\begin{itemize}
\item<1-> Here we present the example from Fukaya-Kato \cite[Example 2.5]{FK} on how we could derive the $p$-adic Kubota-Leopoldt function and Iwasawa Main conjecture from the corresponding big picture above.
\item<2-> Now let $F_n$ be the field $\mathbb{Q}(\zeta_{p^n})^{+}$. Then we use the notation $F_\infty$ to denote the corresponding union throughout all $n$.
\item<3-> Then let $\Lambda$ be $\mathbb{Z}_p[[G]]$, and let $T$ be $\Lambda^\sharp(1)$, where $G$ is the Galois group $\mathrm{Gal}(F_\infty/\mathbb{Q})$.
\item<4-> Now consider $N$ which is largest throughout all the pro-$p$ abelian extension of $F_\infty$ unramified outside all the primes of $F_\infty$ above the prime number $p$.
\item<5-> Then we look at the cohomology groups $H_\mathrm{IW,c}(\mathbb{Z}[1/p],T)$ of the cyclotomic deformation of $\mathbb{Z}_p(1)$ which is also isomorphic to $\varprojlim_{n\rightarrow \infty} HR\Gamma_{c}(\mathcal{O}_{F_n}[1/p],\mathbb{Z}_p(1))$.
\item<6-> By classical Iwasawa theory we have that $H^2_\mathrm{IW,c}(\mathbb{Z}[1/p],T)\overset{\sim}{\rightarrow} \mathrm{Gal}(N/F_\infty)$ and $H^3_\mathrm{IW,c}(\mathbb{Z}[1/p],T)\overset{\sim}{\rightarrow} \mathbb{Z}_p$.
\item<7-> Now we consider the Iwasawa deformation and the corresponding noncommutative Tamagawa number conjecture, so we have an isomorphism:
\begin{align}
\mathrm{Determinant}_\Lambda(0)\overset{\sim}{\longrightarrow}\mathrm{Determinant}_\Lambda R\Gamma_{c}(\mathbb{Z}[1/p],T).	
\end{align}
\item<8-> For $r$ an odd integer which is negative. Then we have then by composing with the $p$-adic period-regulator homomorphism, and by consider the ring homomorphism $\chi_{\mathrm{cyc}}^{1-r}:\Lambda\rightarrow \mathbb{Z}_p$:
\begin{align}
\mathrm{Determinant}_{\mathbb{Q}_p}(0)\overset{\sim}{\longrightarrow}\mathrm{Determinant}_{\mathbb{Q}_p} R\Gamma(\mathbb{Z}[1/p],\mathbb{Q}_p(r))\overset{\sim}{\rightarrow} \mathrm{Determinant}_{\mathbb{Q}_p}(0),
\end{align}
which corresponds to $(1-p^{-r})\zeta(r)$ living in $\mathbb{Q}^\times$, which further corresponds to element in $K_1(\mathbb{Q}_p)$.
\end{itemize}

\noindent{\text{Example: $p$-adic Kubota-Leopoldt and Iwasawa Main Conjecture}}

\begin{itemize}
\item<1-> Now we consider the total fraction $Fr(\Lambda)$ of $\Lambda$ we have then the corresponding isomorphism:
\begin{align}
\mathrm{Determinant}_{Fr(\Lambda)}(0)\overset{\sim}{\longrightarrow}Fr(\Lambda)\otimes\mathrm{Determinant}_\Lambda R\Gamma_{c}(\mathbb{Z}[1/p],T).	
\end{align}
\item<2-> Now it is well-known that  $R\Gamma_{c}(\mathbb{Z}[1/p],T)$
is $\Lambda$-torsion, so we have:
\begin{align}
\mathrm{Determinant}_{Fr(\Lambda)}(0)\overset{\sim}{\longrightarrow}Fr(\Lambda)\otimes\mathrm{Determinant}_\Lambda R\Gamma_{c}(\mathbb{Z}[1/p],T)\overset{\sim}{\longrightarrow} \mathrm{Determinant}_{Fr(\Lambda)}(0).	
\end{align}
\item<3-> This corresponds to an element in $Fr(\Lambda)^\times$ and further an element in the $K$-group $K_1(Fr(\Lambda))$. This is just the corresponding $p$-adic Kubota-Leopoldt element $\zeta_{p,\mathrm{KL}}$.

\item<4-> Then we can consider the corresponding characteristic ideal $(\zeta_{p,\mathrm{KL}})$ generated by the $p$-adic Kubota-Leopoldt element, which is by using the corresponding isomorphism above we have that $(\zeta_{p,\mathrm{KL}})=\mathrm{Characteristic}_\Lambda(R\Gamma_{c}(\mathbb{Z}[1/p],T))$.

\item<5-> Then note that:
\begin{align}
\mathrm{Characteristic}_\Lambda(R\Gamma_{c}(\mathbb{Z}[1/p],T))&\overset{\sim}{\longrightarrow} \mathrm{Characteristic}_\Lambda H^2_\mathrm{IW,c}(\mathbb{Z}[1/p],T) \mathrm{Characteristic}_\Lambda H^3_\mathrm{IW,c}(\mathbb{Z}[1/p],T)^{-1}\\
&\overset{\sim}{\longrightarrow} \mathrm{Characteristic}_\Lambda (\mathrm{Gal}(N/F_\infty) )\mathrm{Characteristic}_\Lambda( \mathbb{Z}_p)^{-1}.	
\end{align}
\item<6-> Then by using the $\zeta$-isomorphism this gives rise to the classical Iwasawa main conjecture \cite{Iwa}:
\begin{align}
(\zeta_{p,\mathrm{KL}})=\mathrm{Characteristic}_\Lambda(G(N/F_\infty))\mathrm{Characteristic}_\Lambda(\mathbb{Z}_p)^{-1}.
\end{align}

\end{itemize}

\noindent{\text{Example: Local $\varepsilon$-Isomorphism}}	

\begin{itemize}
\item<1-> We present the example from Fukaya-Kato \cite[Example 3.6]{FK} in the cyclotomic Iwasawa situation. 
\item<2-> We consider the group $G=\mathrm{Gal}(\mathbb{Q}_p(\zeta_{p^\infty})/\mathbb{Q}_p)$.
\item<3-> Let $\Lambda$ be the Iwasawa algebra $\mathbb{Z}_p[[G]]$.
\item<4-> Let $T$ be the representation of $G_{\mathbb{Q}_p}$ which is given by $\Lambda^{\sharp}(1)$.
\item<5-> Now consider the multiplicative group $\mathbb{Z}_p[\zeta_{p^n}]^\times$, and put $U_n$ to be the corresponding completion by considering the corresponding quotients taking the form of $\mathbb{Z}_p[\zeta_{p^n}]^\times/(\mathbb{Z}_p[\zeta_{p^n}]^\times)^{p^m}$ for $m\geq 0$.
\item<6-> Then one uses the notation $U$ to denote the local units $\varprojlim_{n\rightarrow \infty}U_n$.
\item<7-> Then we have by using the vanishing of Galois cohomology:
\begin{align}
\mathrm{Determinant}_\Lambda(R\Gamma(\mathbb{Q}_p,T))&\overset{\sim}{\longrightarrow}\mathrm{Determinant}_\Lambda(H^1(\mathbb{Q}_p,T))^{-1}\mathrm{Determinant}_\Lambda(H^2(\mathbb{Q}_p,T))\\
&	\overset{\sim}{\longrightarrow}\mathrm{Determinant}_\Lambda(U)^{-1}\mathrm{Determinant}_\Lambda(\mathbb{Z}_p)^{-1}\mathrm{Determinant}_\Lambda(\mathbb{Z}_p)\\
&\overset{\sim}{\longrightarrow}\mathrm{Determinant}_\Lambda(\Lambda)^{-1}\\
&\overset{\sim}{\longrightarrow}\mathrm{Determinant}_\Lambda(T)^{-1}.
\end{align}
\item<8-> Therefore we have $\mathrm{Determinant}_\Lambda(0)\overset{\sim}{\rightarrow} \mathrm{Determinant}_\Lambda(R\Gamma(\mathbb{Q}_p,T))\mathrm{Determinant}_\Lambda(T)$.

\end{itemize}

\newpage

\section{Beyond the \'Etale Picture in the Relative Setting}

\noindent{\text{Nakamura's Local Tamagawa Number Conjecture for $(\varphi,\Gamma)$-Modules}}

\begin{itemize}
\item<1-> Here we present Nakamura's local $\varepsilon$-isomorphism conjecture. We follow his notations here.
\item<2-> $A$ be a reduced affinoid algebra in rigid analytic geometry  with coefficient $\mathbb{Q}_p$.
\item<3-> Now we consider the relative Robba ring with coefficient $A$ which is the ring of Laurent analytic functions defined over the punctured open unit disc relative to $A$, we denote it by $\Pi_{A}$.
\item<4-> Now consider the quasi-coherent sheaves over the relative punctured open unit disc defined above.  
\item<5-> Let $M$ now be the global section of such vector bundle which is a finite projective module over $\Pi_{A}$ which descends to some finite projective module over the global section of some half-open annulus. 
\item<6-> Now we use the same notation $M$ to denote the module carrying further the action of the Frobenius and the quotient $\Gamma$ of the absolute Galois group of $\mathbb{Q}_p$.
\item<7-> In our situation we have the perfect complex $C^*_{\varphi,\Gamma}(M)$ in the derived category of $A$-modules (Kedlaya-Pottharst-Xiao \cite[Theorem, Chapter 1]{KPX}).  
\item<8-> Now since we have the corresponding perfectness so we could define the corresponding determinant in some 'Picard category' (Deligne's categorical group philosophy \cite[Chapitre 4]{De3}). 
	
\item<9-> More or less just some determinant we have the corresponding determinant module $\triangle_A(M)$.
\end{itemize}

\noindent{\text{Local Noncommutative Tamagawa Number Conjecture for $(\varphi,\Gamma)$-Modules}}
\begin{itemize}
\item<1-> Here we need to twist by using some line bundle coming from Kedlaya-Pottharst-Xiao's results \cite[Theorem 6.2.14]{KPX} on the rank 1 Frobenius module.
\item<2-> Nakamura mimicks the Fukaya-Kato's local $\varepsilon$-isomorphism for Galois representations to have the chance to formulate the corresponding $\varepsilon$-isomorphism conjecture for a $(\varphi,\Gamma)$-module $M$. 
\item<3-> (\text{Conjecture, Nakamura, \cite[Conjecture 1.1]{N1}}) There is a unique and canonical way to associate to $M$ an isomorphism (after choosing a basis of $\mathbb{Z}_p(1)$):
\begin{align}
1_A\overset{\sim}{\longrightarrow}	\Delta_A(M)
\end{align}
 which is compatible with extension and functional equation, and compatible (when $A$ is $\mathbb{Q}_p$) with Fukaya-Kato's \'etale $\varepsilon$-isomorphism conjecture, and compatible with the corresponding de-Rham $\varepsilon$-isomorphism conjecture. 
\item<4-> Pal-Z\'abr\'adi \cite{PZ} proposed that one can consider vector bundles over polydiscs and one can consider the corresponding product version of the conjecture and consider the corresponding diagonal embedding to make progress. 
\item<5-> One also has the noncommutative version of the conjecture proposed by Z\"ahringer \cite{Z}.
\end{itemize}

\bibliographystyle{splncs}

\chapter{Survey on $\infty$-Categorical and Homotopical Perspectives}

\section{\text{The Goal}}

In the dissertation \cite{T7} we study the geometric and the representation theoretic aspects of $p$-adic motives, through the study of the corresponding period spectra. The aspects we would like to mention here are as in the following:\\
\begin{itemize}
\justifying
\item<1-> Derived topological de Rham complexes and the derived topological logarithmic de Rham complexes after Bhatt, Gabber, Guo, Illusie, Olsson \cite{B1}, \cite{O1}, \cite{G}, \cite{Ill1}, \cite{Ill2}. Work of Bhatt-Scholze and Li-Liu \cite{BS}, \cite{LL} showed that these are in some sense equivalent to the derived prismatic cohomological complexes. This is kind of crucial since one can construct certain period sheaves in the second item below as in the work of Bhatt, Guo and Li \cite{B1}, \cite{GL}, directly from these complexes. And this is crucial and convenient in the study of singular spaces; \\  

\item<2-> $\mathcal{O}\mathbb{B}_{\mathrm{dR},X}$-sheaves and $\mathcal{O}\mathbb{B}_{\mathrm{dR},\mathrm{log},X}$-sheaves after Diao, Lan, Liu, Scholze, Zhu \cite{DLLZ}, \cite{Scho1}. One should be able to use the first item to construct $\mathcal{O}\mathbb{B}_{\mathrm{dR},\mathrm{log},X}$ for completely intersection rigid analytic space after Bhatt, Guo and Li \cite{B1}, \cite{GL};\\

\item<3-> $\varphi$-$\widetilde{C}_{X}$ sheaves, relative $B$-pairs after Kedlaya-Liu \cite{KL1}, \cite{KL2};\\

\item<4-> Multidimensional Robba rings and sheaves, and multidimensional Frobenius modules and sheaves after Carter-Kedlaya-Z\'abr\'adi and Pal-Z\'abr\'adi \cite{CKZ}, \cite{PZ}. This is towards some progress on $p$-adic Tamagawa number conjecture proposed by Pal-Z\'abr\'adi \cite{PZ};\\

\item<5-> And many other possible motivic period ring spectra such as the derived $I$-adically complete THH and HH of derived $I$-adically complete $\mathbb{E}_1$-rings \cite{NS}.\\

\item<6-> The foundation we rely on could actually be promoted to the level of Bambozzi-Ben-Bassat-Kremnizer, Ben-Bassat-Mukherjee, Bambozzi-Kremnizer, Clausen-Scholze and Kelly-Kremnizer-Mukherjee \cite{BBBK}, \cite{BBM}, \cite{BK}, \cite{CS2}, \cite{CS1}, \cite{KKM}. The project \cite{M} defined derived prismatic cohomology for rings in \cite{CS2}, \cite{CS1}. Also see \cite{T8} and \cite[Chapter 12]{T7} for more discussion around motivic constructions for the rings in \cite{BBBK}, \cite{BBM}, \cite{BK}, \cite{CS2}, \cite{CS1}, \cite{KKM}.\\
\end{itemize}

\newpage
\section{\text{Motivations}}

\indent The objects we studied are mainly crucial in $p$-adic motives, we will consider the analytic geometry of the corresponding spaces over or attached to the period rings mentioned above, with the following key motivations:
\begin{itemize}
\justifying
\item<1-> Noncommutative Tamagawa number conjectures and Noncommutative Iwasawa number conjectures in the sense of Burns, Flach, Fukaya, Kato \cite{BF1}, \cite{BF2}, \cite{FK}; 
\item<2-> Deformation and families of representations of the fundamental groups in analytic and algebraic geometry; 
\item<3-> Deformation and families of generalized motivic structures beyond the representations of fundamental groups;
\item<4-> Noncommutative analytic geometry and deformation theory;
\item<5-> Noncommutative derived analytic geometry and derived deformation theory;
\item<6-> Langlands programs;
\item<7-> Analytic approaches to algebraic topology.
\end{itemize}

\

\noindent{\text{Example of the Rings and Spaces}}\\

\indent Here we discuss a little bit about some examples of the rings or the spaces we are interested:
\begin{itemize}
\justifying
\item<1-> The rigid analytic affinoids or spaces after Tate \cite{Ta}, such as the following ring and the rigid analytic space attached:
\begin{align}
\mathbb{Q}_p\left<T_1,T_2,...,T_d\right>.	
\end{align}
These are actually also the motivation for Huber's adic rings \cite{Hu} which unifies the formal adic rings as in EGA \cite{EGAI1}, \cite{EGAII}, \cite{EGAIII1}, \cite{EGAIII2}, \cite{EGAIV1}, \cite{EGAIV2}, \cite{EGAIV3}, \cite{EGAIV4} and their generic fibers. 
\item<2-> Fukaya-Kato adic ring and generic fibers \cite{FK}, such as\footnote{$Z_1,...,Z_d$ are now free variables, commuting with $T_i,i=1,...,e$.} :
\begin{align}
\mathbb{Z}_p[[T_1,...,T_e]][[Z_1,...,Z_d]], \mathbb{Q}_p\left<T_1,...,T_e\right>\left<Z_1,Z_2,...,Z_d\right>.	
\end{align}
We conjecture there are Fukaya-Kato adic spaces but we have not written down the detail. The $I$-adic completion could be defined on the derived level in this noncommutative situation, namely for $\mathbb{E}_1$-rings\footnote{This is already subtle in the commutative setting as in \cite{BS2}, \cite{Po}, \cite[Section 15.90]{SP} and \cite{Ye}. In the noncommutative setting one also has the chance to take the derived functor of the $J$-adic completion.}.	
\end{itemize}

\begin{itemize}
\justifying
\item<3-> Pseudorigid analytic rings and spaces, for instance the following Tate Huber pair although over $\mathbb{Z}_p$:
\begin{align}
(\mathbb{Z}_p[[t]]\left<\frac{p^a}{t^b}\right>[1/t]\left<T_1,T_2,...,T_d\right>, \mathbb{Z}_p[[t]]\left<\frac{p^a}{t^b}\right>\left<T_1,T_2,...,T_d\right>), (a,b)=1.	
\end{align}

\item<4-> Affinoids over prisms $(A,I)$ of Bhatt-Scholze \cite{BS} (such as $A_\mathrm{inf}(\mathcal{O}_{\mathbb{C}^\flat}),W(\mathbb{F}_p)[[u]],\mathbb{Z}_p[[q-1]]$).

\item<5-> The main difficulty is that in our study of the relative analytic geometry, the base ring is not a field but a very complicated ring. For instance the analytic geometry over the affinoids over Robba rings after Kedlaya-Liu \cite{KL1}, \cite{KL2} is a very complicated subject.  	
\end{itemize}

\

\newpage
\section{\text{$\infty$-Categorical Completion and $\infty$-Categorical Solidification}}

\indent Completions are very important in our study, this is because largely our interested rings carry topology or norms on the rings of definitions which implies that the period rings attached might also carry some topology or norms, therefore we can take the corresponding completion. We actually have some weaker notions on the derived level.

\begin{itemize}
\justifying
\item<1-> As in the work of Bhatt-Scholze \cite{BS} on the corresponding prismatic cohomology, we have the inclusion of the following $\infty$-category from the category of $\mathbb{E}_\infty$-rings attached to some presentable $\infty$-category over some $R$ which are derived $I\subset \pi_0(R)$ to the category of $\mathbb{E}_\infty$-rings attached to some presentable $\infty$-category over some $R$:
\begin{align}
\mathbb{E}_\infty\mathrm{Sp}_{R,I-\mathrm{com}}\rightarrow \mathbb{E}_\infty\mathrm{Sp}_{R}	
\end{align}
with the derived $I$-completion:
\begin{align}
\mathbb{E}_\infty\mathrm{Sp}_{R}	\rightarrow\mathbb{E}_\infty\mathrm{Sp}_{R,I-\mathrm{com}}.
\end{align}
\item<2-> As in the work of Bambozzi-Ben-Bassat-Kremnizer, Ben-Bassat-Mukherjee and Kelly-Kremnizer-Mukherjee \cite{BBBK}, \cite{BBM} and \cite{KKM}, one may have the inclusion of the following $\infty$-categories from the category of $\mathbb{E}_\infty$-rings attached to the simplicial ind-Banach sets over some Banach $R$ to the category of $\mathbb{E}_\infty$-rings attached to the simplicial ind-normed sets over some $R$:
\begin{align}
\mathbb{E}_\infty\mathrm{Simp}_{R,\mathrm{ind-Ban}}\rightarrow \mathbb{E}_\infty\mathrm{Simp}_{R,\mathrm{ind-Nor}}	
\end{align}
which admits a left adjoint functor namely the derived Banach completion.

\end{itemize}

\begin{itemize}
\justifying
\item<3-> From Clausen-Scholze \cite{CS1}, \cite{CS2} we have the solidification functor:
\begin{align}
(.)^{\mathbb{L}\mathrm{Solidified}}: D_{\mathrm{condensed, ab}}\rightarrow D_{\mathrm{condensed, ab, solid}}	
\end{align}
which is the left adjoint of the corresponding inclusion:
\begin{align}
D_{\mathrm{condensed, ab, solid}}\rightarrow D_{\mathrm{condensed, ab}}.	
\end{align}
Solids are characterized as the the category closed under the limits, colimits and extensions, which is somewhat abstract characterization of complete objects.

\item \mbox{\text{($\infty$-Categorical Completion and $\infty$-Categorical Solidification)}} In fact one will have the chance to take many different types of completions, for instance the Hodge-completion or parallely the Nygaard-completion is happening in some filtered derived $\infty$-category which is very significant in the development of derived prismatic cohomology \cite{BS}. Also we have the following foundations. First we have from \cite[Proposition 1.6]{BBBK} the Banach completion on the Banach sets over some Banach ring $R$ which is the left adjoint functor of the inclusion:
\begin{align}
\mathrm{BanSets}_R\rightarrow \mathrm{NormSets}_R.
\end{align}
One may also consider the corresponding functor on the derived level between simplicial ind-Banach and ind-Normed sets:
\begin{align}
\mathrm{Simp}\mathrm{Ind}\mathrm{BanSets}_R\rightarrow \mathrm{Simp}\mathrm{Ind}\mathrm{NormSets}_R,
\end{align}
which may further be extended to the simplicial ind-Banach and ind-Normed abelian groups, rings and modules. And we have on the derived level the corresponding derived $I$-completion over some $\mathbb{E}_\infty$ ring $R$ where $I\subset \pi_0(R)$. When $I$ is finitely generated, if this happens in the derived $\infty$-category $\mathbb{D}_R$ as in \cite[Chapter 1 Notation]{BS} and \cite[Chapter 3]{BS2}, then this for instance could be regarded as the left adjoint functor of the inclusion:
\begin{align}
h\mathbb{D}_{R,I-\mathrm{com}}\rightarrow h\mathbb{D}_R,
\end{align}
when $R$ is classical as in \cite[Tag 091N, Section 15.90]{SP}. Moreover we have the derived $I$-completion on the $\infty$-categorical level to get derived $I$-complete objects:
\begin{align}
\mathbb{D}_R\rightarrow  \mathbb{D}_{R,I-\mathrm{com}}. 
\end{align}
Also we have from \cite[Theorem 5.8]{CS1} the derived solidification functor $(.)^{\mathbb{L}\mathrm{Solidified}}$ on the derived level which is the left adjoint of the following inclusion from derived category of solids to the derived category of condensed abelian groups:
\begin{align}
D_{\mathrm{condensed},\mathrm{ab},\mathrm{solid}}\rightarrow D_{\mathrm{condensed},\mathrm{ab}}.
\end{align}
And furthermore one can consider the corresponding functors on the derived $\infty$-categorical level (also see \cite[Chapter 1]{Lu1}):
\begin{align}
\mathbb{D}_{\mathrm{condensed},\mathrm{ab},\mathrm{solid}}\rightarrow \mathbb{D}_{\mathrm{condensed},\mathrm{ab}}.
\end{align}
Solids are abstractly characterized in \cite{CS1} as those condensed abelian groups closed under limits, colimits and extensions. Therefore solidification is abstract characterization of completion. Solids may be very important in Iwasawa theory when we want to take solidified tensor products under some absolute topology. For instance we may take the product $\mathcal{F}\otimes^{\mathbb{L}\mathrm{Solidified}}P[1/I]^\wedge_p$ of some Fr\'echet algebra $\mathcal{F}$ such as in \cite{ST} with certain ring $P[1/I]^\wedge_p$ (or any similar one) produced from some prism $(P,I)$ in \cite{BS} under this framework. At least for instance in deformation theory happening over the Robba ring $\widetilde{\Pi}^I_{R}\otimes^{\mathbb{L}\mathrm{Solidified}}\mathcal{F}$ after Kedlaya-Liu \cite{KL1} and \cite{KL2}, sheafiness has been not a problem anymore after \cite{BK} and \cite{CS2}, which implies that one should look at ind-Banach sets instead of Banach sets and pro-\'etale solids instead of complete sets.

\end{itemize}

\

\noindent{\text{$\infty$-Categorical Iwasawa-Prismatic Theory with $\infty$-Categorical Hodge-Iwasawa Theory}}\\

\indent We now consider the analog of Hodge-Iwasawa theory, which will represent the main idea of our consideration in Hodge-Iwasawa theory \cite{T1}, \cite{T2}, for more discussion see \cite[Chapter 1]{T7}. We consider the following discussion after \cite{BS}, \cite{BS3} and \cite{3Sch2}:

\begin{itemize}
\justifying
\item<1-> Let $A$ be a quasisyntomic formal ring over $\mathbb{Z}_p$, for instance the ring of integer of a local field. Let $X$ be $\mathrm{Spf}A$.
\item<2-> We have the prismatic site $(X_{\mathrm{prim}},\mathcal{O}_{X_{\mathrm{prim}}})$. 
\item<3-> We have the quasisyntomic site $(X_{\mathrm{qsyn}},\mathcal{O}_{X_{\mathrm{qsyn}}})$, recall for any quasiregular semiperfectoid $A$ in the site, we consider the associated prism $(P_*,I_{P_*})(A)$, then we have $\mathcal{O}_{X_{\mathrm{qsyn}}}(A):=P_*(A)$.
\item<4-> Now we consider $A=\mathcal{O}_K$ of some $p$-adic local field. Work of Bhatt-Scholze \cite{BS}, \cite{BS3} and \cite{3Sch2} tells us that the category of all the finite locally free $F$-crystals over the prismatic site $(X_{\mathrm{prim}},\mathcal{O}_{X_{\mathrm{prim}}}[1/I_{\mathcal{O}_{X_{\mathrm{prim}}}}]^\wedge_p)$ is equivalent to the category of all the $\mathbb{Z}_p$-representations of $G_K$, while the category of all the finite locally free $F$-crystals over the prismatic site $(X_{\mathrm{prim}},\mathcal{O}_{X_{\mathrm{prim}}}[1/I_{\mathcal{O}_{X_{\mathrm{prim}}}}]^\wedge_p[1/p])$ is equivalent to the category of all the $\mathbb{Q}_p$-representations of $G_K$.

\end{itemize}

\begin{itemize}
\justifying
\item<5-> Now we take an Iwasawa algebra $\mathbb{Z}_p[[G]]$ of some compact $p$-adic Lie group $G$, and consider the sites:
\begin{align}
(X_{\mathrm{prim}},\mathcal{O}_{X_{\mathrm{prim}}}[1/I_{\mathcal{O}_{X_{\mathrm{prim}}}}]^\wedge_p{\otimes}^{\mathbb{L}\mathrm{com}}\mathbb{Z}_p[[G]]), (X_{\mathrm{qsyn}},\mathcal{O}_{X_{\mathrm{qsyn}}}[1/I_{\mathcal{O}_{X_{\mathrm{qsyn}}}}]^\wedge_p	{\otimes}^{\mathbb{L}\mathrm{com}}\mathbb{Z}_p[[G]]),\\
(X_{\mathrm{prim}},\mathcal{O}_{X_{\mathrm{prim}}}[1/I_{\mathcal{O}_{X_{\mathrm{prim}}}}]^\wedge_p{\otimes}^{\mathbb{L}\mathrm{com}}\mathbb{Z}_p[[G]][1/p]), (X_{\mathrm{qsyn}},\mathcal{O}_{X_{\mathrm{qsyn}}}[1/I_{\mathcal{O}_{X_{\mathrm{qsyn}}}}]^\wedge_p	{\otimes}^{\mathbb{L}\mathrm{com}}\mathbb{Z}_p[[G]][1/p]),
\end{align}
after taking any possible completion in the derived sense.
\item<6-> Then this kind of deformation will be the main idea of the so-called Iwasawa deformation. Clausen-Scholze's solidified tensor product \cite{CS1}, \cite{CS2} may be relevant here when we want to take the tensor product $\mathcal{F}{\otimes^{\mathbb{L}\mathrm{Solidified}}}\mathcal{O}_P[1/I_P]^\wedge_p$ of an ind-Fr\'echet algebra $\mathcal{F}$ with any sheaf $\mathcal{O}_P[1/I_P]^\wedge_p$ as above constructed from some sheaf of prism $(\mathcal{O}_P,\mathcal{I}_P)$, or when we want to take the tensor product of Kedlaya-Liu Robba ring $\widetilde{\Pi}_R^I$ with solid condensed rings such as any adic ring.	
\end{itemize}

\newpage

\section{First Scope of the Study: $\infty$-Categorical Algebraic Hodge-Iwasawa Sheaves}

\indent The first scope of the study is around the corresponding rigid families of the corresponding Frobenius modules and the quasicoherent sheaves over the corresponding schematic Fargues-Fontaine curves in the relative $p$-adic Hodge theory of Kedlaya-Liu and Kedlaya-Pottharst \cite{KL1}, \cite{KL2}, \cite{KP}, directly aiming at the corresponding deformations of profinite fundamental groups. 

\begin{itemize}
\justifying
\item<1-> Directly generalizing work of Kedlaya-Liu and Kedlaya-Pottharst \cite{KL1}, \cite{KL2}, \cite{KP} we have the following theorem:

\item<2-> (\text{Theorem, \cite[Theorem 1.3, Theorem 1.4]{T1}, \cite[Proposition 4.25, Proposition 4.26]{T2}})
Let $(R,R^+)$ be a char $p>0$ perfect uniform adic Banach ring which is Tate and $A$ is a rigid analytic affinoid over some analytic field ($\mathbb{Q}_p$ or $\mathbb{F}_p((t))$). Then we have:\\
1. There is an equivalence between the category of vector bundles over the schematic Fargues-Fontaine curve $\mathrm{Proj}_{\widetilde{\Pi}_{R,A}}$, the category of Frobenius finite projective modules over $\widetilde{\Pi}_{R,A}$, the category of Frobenius finite projective bundles over $\widetilde{\Pi}_{R,A}$ as well as the category of Frobenius finite projective modules over $\widetilde{\Pi}^\infty_{R,A}$;\\
2. There is an equivalence between the category of pseudocoherent sheaves over the schematic Fargues-Fontaine curve $\mathrm{Proj}_{\widetilde{\Pi}_{R,A}}$ and the category of algebraic Frobenius pseudocoherent modules over $\widetilde{\Pi}_{R,A}$.\\

\item<3-> The algebraicity here means that we do not know if the natural topology coming from the finite free covering will be maintained to be complete as well, certainly for the theory of the pseudocoherent sheaves over the schemes above is also purely algebraic.

\end{itemize}

\newpage

\section{Second Scope of the Study: $\infty$-Categorical Analytic Hodge-Iwasawa Sheaves}
	
\indent The second scope of the study is around the corresponding rigid families of the corresponding Frobenius modules and the pseudocoherent sheaves over the corresponding adic Fargues-Fontaine curves in the relative $p$-adic Hodge theory of Kedlaya-Liu and Kedlaya-Pottharst \cite{KL1}, \cite{KL2}, \cite{KP}. 
\begin{itemize}
\justifying
\item<1-> Directly generalizing work of Kedlaya-Liu and Kedlaya-Pottharst \cite{KL1}, \cite{KL2}, \cite{KP} we have the following theorem:

\item<2-> (\text{Theorem, \cite[Proposition 5.44]{T2}})
Let $(H_\bullet,H_\bullet^+)$ be a noetherian, weakly-decompleting, decompleting, finite \'etale, perfectoid tower from Kedlaya-Liu. Assume that this is Galois with Galois group $\Gamma$. Let $A$ be a sousperfectoid rigid analytic affinoid. Then we have the categories of $(\varphi,\Gamma)$ finite projective modules over $\widetilde{\Pi}_{(H_\bullet,H_\bullet^+),A}$, $\breve{\Pi}_{(H_\bullet,H_\bullet^+),A}$, ${\Pi}_{(H_\bullet,H_\bullet^+),A}$, the the categories of $(\varphi,\Gamma)$ finite projective budles over $\widetilde{\Pi}_{(H_\bullet,H_\bullet^+),A}$, $\breve{\Pi}_{(H_\bullet,H_\bullet^+),A}$, ${\Pi}_{(H_\bullet,H_\bullet^+),A}$, the category of $\Gamma$-equivariant vector bundles over adic FF-curve and the category of all the Frobenius sheaves over the corresponding Robba sheaves in the pro-\'etale sites, are all equivalent to each other.

\item<3-> Here the product happens between the Kedlaya-Liu Robba ring with $A$ is the Banach completed tensor product, one may want to consider the theory by taking the product $\otimes^{\mathbb{L}\mathrm{Solidified}}$. This may be very crucial when we want to take the product with any pseudorigid analytic affinoid as in Bellovin's work to treat some general eigenvarieties \cite{3Bel1}, \cite{Bel2}. Therefore this will be helpful when we want to take the product of rings having different topologies.

\end{itemize}


\begin{itemize}
\justifying
\item<4-> (\text{Conjecture})
Let $(H_\bullet,H_\bullet^+)$ be a noetherian, weakly-decompleting, decompleting, finite \'etale, perfectoid tower from Kedlaya-Liu. Assume that this is Galois with Galois group $\Gamma$. Let $A$ be a simplicial topological algebra\footnote{For instance one can consider pseudorigid affinoids and more general topological rings carrying linear topology.}. Then we have the categories of $(\varphi,\Gamma)$ finite projective modules over $\widetilde{\Pi}_{(H_\bullet,H_\bullet^+)}\otimes^{\mathbb{L}\mathrm{Solidified}}A$, $\breve{\Pi}_{(H_\bullet,H_\bullet^+)}\otimes^{\mathbb{L}\mathrm{Solidified}}A$, ${\Pi}_{(H_\bullet,H_\bullet^+)}\otimes^{\mathbb{L}\mathrm{Solidified}}A$, the the categories of $(\varphi,\Gamma)$ finite projective budles over $\widetilde{\Pi}_{(H_\bullet,H_\bullet^+)}\otimes^{\mathbb{L}\mathrm{Solidified}}A$, $\breve{\Pi}_{(H_\bullet,H_\bullet^+)}\otimes^{\mathbb{L}\mathrm{Solidified}}A$, ${\Pi}_{(H_\bullet,H_\bullet^+)}\otimes^{\mathbb{L}\mathrm{Solidified}}A$, are all equivalent to each other.	
\item<5-> Directly generalizing work of Kedlaya-Liu and Kedlaya-Pottharst \cite{KL1}, \cite{KL2}, \cite{KP} we have the following theorem:

\item<6-> (\text{Theorem, \cite[Proposition 5.51]{T2}})
Let $(H_\bullet,H_\bullet^+)$ be a noetherian, weakly-decompleting, decompleting, finite \'etale, perfectoid tower from Kedlaya-Liu. Assume that this is Galois with Galois group $\Gamma$, and we assume the tower is $F$-projective. Let $A$ be a sousperfectoid rigid analytic affinoid. Then we have:\\
I. The categories of $(\varphi,\Gamma)$ pseudocoherent modules over $\widetilde{\Pi}_{(H_\bullet,H_\bullet^+),A}$, $\breve{\Pi}_{(H_\bullet,H_\bullet^+),A}$, ${\Pi}_{(H_\bullet,H_\bullet^+),A}$;\\ 
II. The the category of $(\varphi,\Gamma)$ pseudocoherent bundles over $\widetilde{\Pi}_{(H_\bullet,H_\bullet^+),A}$, $\breve{\Pi}_{(H_\bullet,H_\bullet^+),A}$, ${\Pi}_{(H_\bullet,H_\bullet^+),A}$;\\ 
III. The category of $\Gamma$-equivariant pseudocoherent sheaves over the adic FF-curve;\\
IV. The categories of all the Frobenius sheaves over the corresponding Robba sheaves in the pro-\'etale sites, are all equivalent to each other.

\end{itemize}


\begin{itemize}
\justifying
\item<7-> The conjectures for any simplicial Banach deformation $A$ and the rings $\widetilde{\Pi}_{(H_\bullet,H_\bullet^+)}\otimes^{\mathbb{L}\mathrm{Solidified}} A$, $\breve{\Pi}_{(H_\bullet,H_\bullet^+)}\otimes^{\mathbb{L}\mathrm{Solidified}} A$, ${\Pi}_{(H_\bullet,H_\bullet^+)}\otimes^{\mathbb{L}\mathrm{Solidified}} A$ could be formulated and possibly proved by directly applying Clausen-Scholze construction and descent \cite{CS1}, \cite{CS2}.  

\item<8-> This would be already very interesting in the situation where we have that the tower is over a local field, which will be directly related to Galois deformation theory. In this situation we have a well-defined $\infty$-derived category of derived ind-coherent $(\varphi,\Gamma)$-sheaves in Lurie's book on higher algebra \cite{Lu1}. This is because we have that the ind-category of coherent sheaves in our currently situation is Grothendieck abelian category admitting enough injective objects.

\item<9-> We are searching derived $\infty$-categories as in the work of Deligne, Fukaya-Kato, Muro-Tonks, Witte \cite{De3}, \cite{FK}, \cite{MT}, \cite{Wit1}. Clausen-Scholze derived category $D(\mathrm{Mod}_{\mathrm{Solidified},\Pi})$ and Fargues-Scholze category $D_\mathrm{Solidified}(*,\Pi)$ might be very crucial here \cite{CS1}, \cite{CS2}, \cite{FS}. The latter is for small $v$-stack, but we can consider the seminormal rigid analytic spaces which could be regarded as small $v$-stack, so the latter is still applicable. And we have the categories $D_\text{\'et}(*,\Pi)$ and $D_\text{qpro\'et}(*,\Pi)$ in Scholze's work on \'etale cohomology of small $v$-stacks \cite{3Sch3}. Note that they are not constructed directly at least on the level of small $v$-stacks, instead they are well-defined subcategories of $D(*_v,\Pi)$ which is definitely site theoretic straightforward construction.

\end{itemize}

\newpage

\section{Third Scope of the Study: $\infty$-Categorical Generalized Hodge-Iwasawa Sheaves}

\indent With the motivation from Pal-Z\'abr\'adi program \cite{PZ} partially aiming at making progress towards the local $p$-adic Tamagawa Number conjecture after Nakamura \cite{N2} on the level of not necessarily \'etale Frobenius modules carrying some Lie group actions, we study some generalized version of the structures after Kedlaya-Pottharst-Xiao \cite{KPX} to allow multidimensional Robba rings and multidimensional Frobenius and Lie group structures. 
\begin{itemize}
\justifying
\item<1-> We consider a finite multi radii $r_I\subset [0,1)^{|I|}$ with $I$ a finite set, then consider the corresponding ring of analytic functions with coefficients in some rigid affinoid $A$ over $\mathbb{Q}_p$ over the rigid multi annulus $[r_1,1)\times...\times [r_{|I|},1)$, which we denote it by $\Pi_{\mathrm{an},[r_1,1)\times...\times [r_{|I|},1),I,A}$. 

\item<2-> Then we take the corresponding injective limit throughout all the multi radii as such to achieve the full multidimensional Robba ring with several variables in our situation:
\begin{align}
\Pi_{\mathrm{an},\mathrm{con},I,A}:=\varinjlim \Pi_{\mathrm{an},[r_1,1)\times...\times [r_{|I|},1),I,A}.	
\end{align}

\item<3-> Relative to $A$ this kind of rings carry the multidimensional Frobenius actions and multidimensional $\Gamma$-action, as well as the multidimensional 'inverse-Frobenius' actions. Applying Kiehl's result \cite[Satz 5.2]{Kie} on strictly completely continuous mappings extended by Kedlaya-Liu \cite{KL3}, we have the following theorem.
	
\end{itemize}

\begin{itemize}

\justifying
\item<4-> (\text{Theorem, \cite[Theorem 1.5]{T3}}) The $(\varphi_I,\Gamma_I)$-complex of any $(\varphi_I,\Gamma_I)$-finite projective module over $\Pi_{\mathrm{an},\mathrm{con},I,A}$ (which is required to be a vector bundle over some $\Pi_{\mathrm{an},[r_1,1)\times...\times [r_{|I|},1),I,A}$) lives in $\mathbb{D}^\flat_\mathrm{perf}(A)$. The Iwasawa complex namely the $\varphi^{-1}_I$-complex of any $(\varphi_I,\Gamma_I)$-finite projective module over $\Pi_{\mathrm{an},\mathrm{con},I,\mathbb{Q}_p}$ (which is required to be a vector bundle over some $\Pi_{\mathrm{an},[r_1,1)\times...\times [r_{|I|},1),I,\mathbb{Q}_p}$) lives in $\mathbb{D}^\flat_\mathrm{perf}(\Pi_{\mathrm{an},[0,1)\times...\times [0,1),I,\mathbb{Q}_p}(\Gamma_I))$.

\item<5-> The concrete property of Iwasawa cohomology groups needs more work, actually we only have that they are coherent sheaves over some non-compact rigid space with $A=\mathbb{Q}_p$. By extending this, one needs to modify some argument from Z\'abr\'adi \cite[Theorem 3.10]{Z} which is conjectured by Kedlaya to be true in even relative situation.\\	
\end{itemize}

\newpage
\section{Fourth Scope of the Study: Sheafiness and $\infty$-Analytic Stacks}

Sheafiness is very crucial in our study above, since the structure presheaves have to be sheaves in order to perform the descent by using results of Kedlaya-Liu \cite[Theorem 2.7.7]{KL1}, \cite[Theorem 2.5.5, Theorem 2.5.14, Section 2.6, Theorem 3.4.8]{KL2}. 

\begin{itemize}
\justifying
\item<1-> One can actually resolve this by using the foundations from Bambozzi-Ben-Bassat-Kremnizer, Ben-Bassat-Mukherjee, Bambozzi-Kremnizer, Clausen-Scholze and Kelly-Kremnizer-Mukherjee \cite{BBBK}, \cite{BBM}, \cite{BK}, \cite{CS2}, \cite{CS1}, \cite{KKM}. Namely one looks at the following inclusions:
\begin{align}
\mathrm{PreAdic}\hookrightarrow \infty-\mathrm{AnStack}_{\mathrm{BBBK},\infty-\mathrm{ringed}},\\
\mathrm{PreAdic}\hookrightarrow \infty-\mathrm{AnStack}_{\mathrm{CS},\infty-\mathrm{ringed}}.
\end{align}

\item<2-> The stacks are $\infty$-ringed which are allowing us to study sheafiness up to higher homotopy, which is the key point to enlarge the categories of analytic spaces.
 
\item<3-> Also one can do this directly through some argument due to Kedlaya on the corresponding lifting splittings through the involved descent morphism. 

\item<4-> Namely whenever we have a corresponding strict exact sequence of Banach rings giving the glueing square where $R_2\rightarrow R_{12}$ has dense image:
\[
\xymatrix@R+0pc@C+0pc{
0  \ar[r]\ar[r]\ar[r] &R \ar[r]\ar[r]\ar[r] &R_1\oplus R_2 \ar[r]\ar[r]\ar[r] &R_{12}  \ar[r]\ar[r]\ar[r] & 0.
}
\] 
Then we have given the finite projective descent datum over $R_1,R_2,R_{12}$, then the equalizer over $R$ is finite projective with desired base change property over $R_1$ and $R_2$. We have the following proofs.

\end{itemize}

\begin{itemize}
\justifying
\item<5-> This is based on the foundation from Kedlaya-Liu and relied on some argument due to Kedlaya, where one lifts the corresponding splittings happening over $R_1$ and $R_2$ to the one over $R$.
\item<6-> One just takes the embedding of the category of adic spaces to the category of Clausen-Scholze, then arrange the modules as sheaves over the corresponding sheaves of animated rings. Then the result follows directly.
\item<7-> This could be extended to the corresponding noncommutative situation since we are considering the deformation through some noncommutative quotient of the free nc Tate algebra $\mathbb{Q}_p\left<Z_1,...,Z_d\right>$\footnote{\justifying In fact we will believe that one can define some $\infty$-noncommutative analytic stacks out of the category $\mathrm{Simp}\mathrm{Ind}\mathrm{BanSets}$ of Bambozzi-Ben-Bassat-Kremnizer and the category $D(\mathrm{Solids})$ of Clausen-Scholze, in some sense parallel to Kontsevich-Rosenberg \cite{3KR1}, \cite{3KR2}. Bambozzi-Ben-Bassat-Kremnizer and Clausen-Scholze stacks are fibered over $\infty$-commutative rings.}. 
\item<8-> Furthermore one can extend this to the derived deformation.  

\item<9-> (\text{Example})
This is the derived rational localization due to Bambozzi-Kremnizer which also produces a corresponding animated ring of Clausen-Scholze, we take a commutative Banach ring over $\mathbb{Q}_p$, then for any $f_1,...,f_d,g$ generating the unit ideal, we take the algebraic Koszul complex attached to these elements for $R$, namely we take $R^h_{f_1,...,f_d,g}:=R\left<T_1,...,T_d\right>/^\mathbb{L}(f_1-gT_1,...,f_d-gT_d)$. Then we take the corresponding complete tensor product with the Robba rings above we have $\widetilde{\Pi}^I_{R}\widehat{\otimes}_{\mathbb{Q}_p}^\mathbb{L}R\left<T_1,...,T_d\right>/^\mathbb{L}(f_1-gT_1,...,f_d-gT_d)$ for some closed interval. This may be useful in derived Galois deformation theory after Galatius-Venkatesh \cite{GV}.

\end{itemize}

\begin{itemize}
\justifying
\item<1-> Then we look at the corresponding glueing sequence as above in more general form with overlapped intervals $I_1,I_2$ with union $I$:
\[
\xymatrix@R+0pc@C+0pc{
\widetilde{\Pi}^I_{R,R^h_{f_1,...,f_d,g}} \ar[r]\ar[r]\ar[r] &\widetilde{\Pi}^{I_1}_{R,R^h_{f_1,...,f_d,g}}\oplus \widetilde{\Pi}^{I_2}_{R,R^h_{f_1,...,f_d,g}} \ar[r]\ar[r]\ar[r] &\widetilde{\Pi}^{I_1\cap I_2}_{R,R^h_{f_1,...,f_d,g}},
}
\] 
which produces the corresponding desired strictly exact sequences for the homotopy groups:
\[
\xymatrix@R+0pc@C-0.2pc{
0\ar[r]\ar[r]\ar[r]  &\pi_k\widetilde{\Pi}^I_{R,R^h_{f_1,...,f_d,g}} \ar[r]\ar[r]\ar[r] &\pi_k\widetilde{\Pi}^{I_1}_{R,R^h_{f_1,...,f_d,g}}\oplus \pi_k\widetilde{\Pi}^{I_2}_{R,R^h_{f_1,...,f_d,g}} \ar[r]\ar[r]\ar[r] &\pi_k\widetilde{\Pi}^{I_1\cap I_2}_{R,R^h_{f_1,...,f_d,g}} \ar[r]\ar[r]\ar[r]  &0,
}
\]
where $k\in \mathbb{Z}_{\geq 0}$. 

\item<2-> Then one can glue the corresponding finite projective objects along these. Again here one can apply Kedlaya's argument on lifting specific splittings, or one can just apply geometric result from Clausen-Scholze descent \cite[Theorem 14.9, Remark 14.10]{CS2} where one definitely has the chance to get also the descent over $\widetilde{\Pi}^I_{R}{\otimes}_{\mathbb{Q}_p}^{\mathbb{L}\mathrm{Solidified}}G$ in some parallel way with considering any adic ring $G$. Certainly it does not matter if the rings have higher homotopies or not.

\end{itemize}

\newpage

\section{Fifth Scope of the Study: Deformation of Functional Analytic Sheaves}

\indent The situation in noncommutative setting is actually very complicated due to the fact that we do not have desired analytic toposes in hand at least at this moment. 
\begin{itemize}

\justifying

\item<1-> But we can deform the corresponding structure sheaves (in the sheafy situation) over analytic field by a noncommutative Banach ring over the same analytic field. 

\item<2-> Consider $X=\mathrm{Spa}(R,R^+)$ and adic affinoid space over some analytic field $F$, which is not just a pre-adic one. 

\item<3-> Then we deform the sheaf by a noncommutative Banach ring $Z$ over $F$ by taking the complete tensor product namely $\mathcal{O}_X\widehat{\otimes}Z$. Then we have after Kedlaya-Liu with the same method \cite[Theorem 2.5.5]{KL2}:

\item<4-> (\text{Theorem, \cite{T4}, \cite[Theorem 1.1]{T5}, \cite{T6}})
The global section functor realizes an equivalence between the category of all the $Z$-equivariant stably pseudocoherent modules over $R\widehat{\otimes}_F Z$ and the category of all the $Z$-equivariant stably pseudocoherent sheaves over $\mathcal{O}_X\widehat{\otimes}Z$.

\item<5-> We can extend this to the \'etale, pro-\'etale sites and quasi-Stein situation as well \cite[Theorem 2.5.14, Section 2.6, Theorem 3.4.8]{KL2}, \cite{T4}, \cite{T5}, \cite{T6}.

\end{itemize}

\newpage

\section*{Acknowledgements}
This survey article is based on the author's dissertation \cite{T7} and beyond. We would like to thank Professor Kedlaya for the suggestions and discussion around topics discussed here.

\bibliographystyle{splncs}

\end{document}